\documentclass[11pt,dvips]{article}
\usepackage{amsmath,amsfonts,amssymb,amsthm,graphicx,verbatim,subfig}
\usepackage[all]{xy}
\usepackage[a4paper,left=25mm]{geometry}
\ifx\pdftexversion\undefined

\usepackage[a4paper,colorlinks,link
=black,filecolor=black,citecolor=black,urlcolor=black,pdfstartview=FitH]{hyperref}
\else

\usepackage[a4paper,colorlinks,linkcolor=black,filecolor=black,citecolor=black,urlcolor=black,pdfstartview=FitH]{hyperref}
\fi

\font\sixbb=msbm6
\font\eightbb=msbm8
\font\twelvebb=msbm10 scaled 1095
\newfam\bbfam
\textfont\bbfam=\twelvebb \scriptfont\bbfam=\eightbb
                           \scriptscriptfont\bbfam=\sixbb
\def\bb{\fam\bbfam\twelvebb}
\newcommand{\Rea}{{\bb R}}
\newcommand{\Com}{{\bb C}}
\newcommand{\Int}{{\bb Z}}

\newcommand{\FF}{{\bb F}}

\newcommand{\PP}{{\bb P}}

\newtheorem{theorem}{\bf Theorem}[section]
\newtheorem{claim}[theorem]{\bf Claim}

\newtheorem{proposition}[theorem]{\bf Proposition}

\newcommand{\enp}{\begin{flushright} $\Box$ \end{flushright}}
\newcommand{\beq}[0]{\begin{equation}}
\newcommand{\enq}[0]{\end{equation}}

\newcommand{\dn}{\Delta_{n-1}}

\newcommand{\thh}{\tilde{H}}
\newcommand{\supp}{{\rm supp}}

\newcommand{\pr}{{\rm Pr}}

\newcommand{\phat}{\widehat{\phi}}

\newcommand{\psihat}{\widehat{\psi}}
\newcommand{\phihat}{\widehat{\phi}}
\newcommand{\img}{{\rm Im} \,}
\newcommand{\ccl}{\mathcal{L}}
\newcommand{\tcl}{\widetilde{\ccl}}

\newcommand{\ghat}{\widehat{G}}
\newcommand{\ghatk}{\ghat^k}

\title{Spectral Expansion of Random Sum Complexes}
\begin{document}
\author{Orr Beit-Aharon \thanks{College of Computer and Information Science,
Northeastern University,
Boston MA 02115, USA. e-mail:
o.beitaharon@northeastern.edu~. } \and Roy Meshulam\thanks{Department of Mathematics,
Technion, Haifa 32000, Israel. e-mail:
meshulam@math.technion.ac.il~. Supported by ISF and GIF grants.} }
\maketitle
\pagestyle{plain}
\begin{abstract}
Let $G$ be a finite abelian group of order $n$ and let $\dn$ denote the $(n-1)$-simplex on the vertex set $G$.
The \emph{sum complex} $X_{A,k}$ associated to a subset $A \subset G$ and $k < n$, is the $k$-dimensional simplicial complex obtained
by taking the full $(k-1)$-skeleton of $\dn$ together with all $(k+1)$-subsets $\sigma \subset G$ that satisfy $\sum_{x \in \sigma} x \in A$. Let $C^{k-1}(X_{A,k})$ denote the space of complex valued $(k-1)$-cochains of $X_{A,k}$.
Let $L_{k-1}:C^{k-1}(X_{A,k}) \rightarrow C^{k-1}(X_{A,k})$ denote the reduced $(k-1)$-th Laplacian
of $X_{A,k}$, and let $\mu_{k-1}(X_{A,k})$ be the minimal eigenvalue of $L_{k-1}$.
\\
It is shown that if $k \geq 1$ and $\epsilon>0$ are fixed, and $A$ is a random subset of $G$ of size $m=\lceil\frac{4k^2\log n}{\epsilon^2}\rceil$, then
\[
\pr\big[~\mu_{k-1}(X_{A,k}) < (1-\epsilon)m~\big] =O\left(\frac{1}{n}\right).
\]
\end{abstract}

\section{Introduction}
\label{s:intro}
\begin{figure}
  \subfloat[A $7$-point triangulation of $\Rea\PP^2$]
  {\label{fig:rp2}
  \scalebox{0.6}{\includegraphics{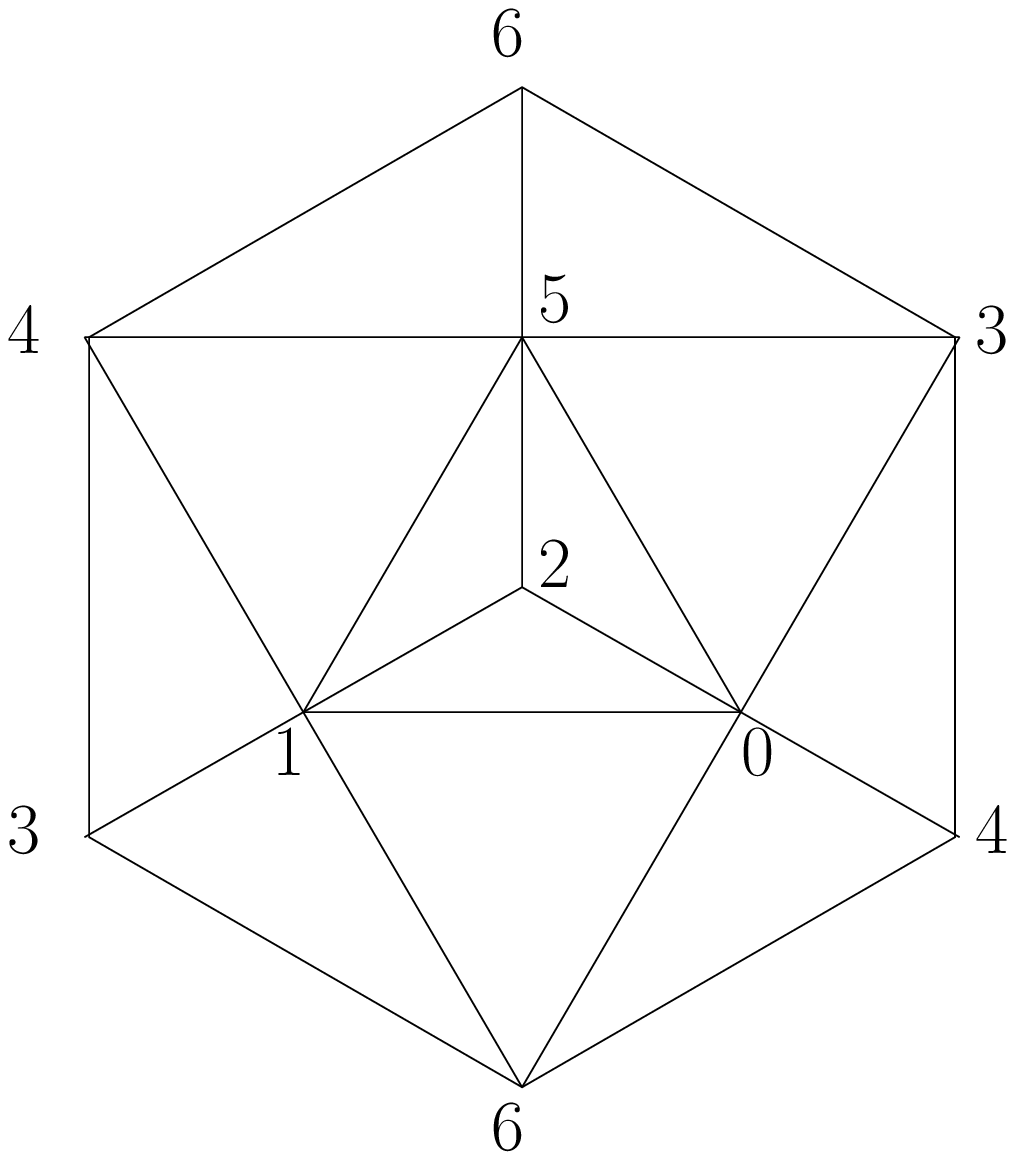}}}
  \hspace{30pt}
  \subfloat[$X_{A,2}$ for $A=\{0,1,3\} \subset \Int_7$]
  {\label{fig:xa}
  \scalebox{0.6}{\includegraphics{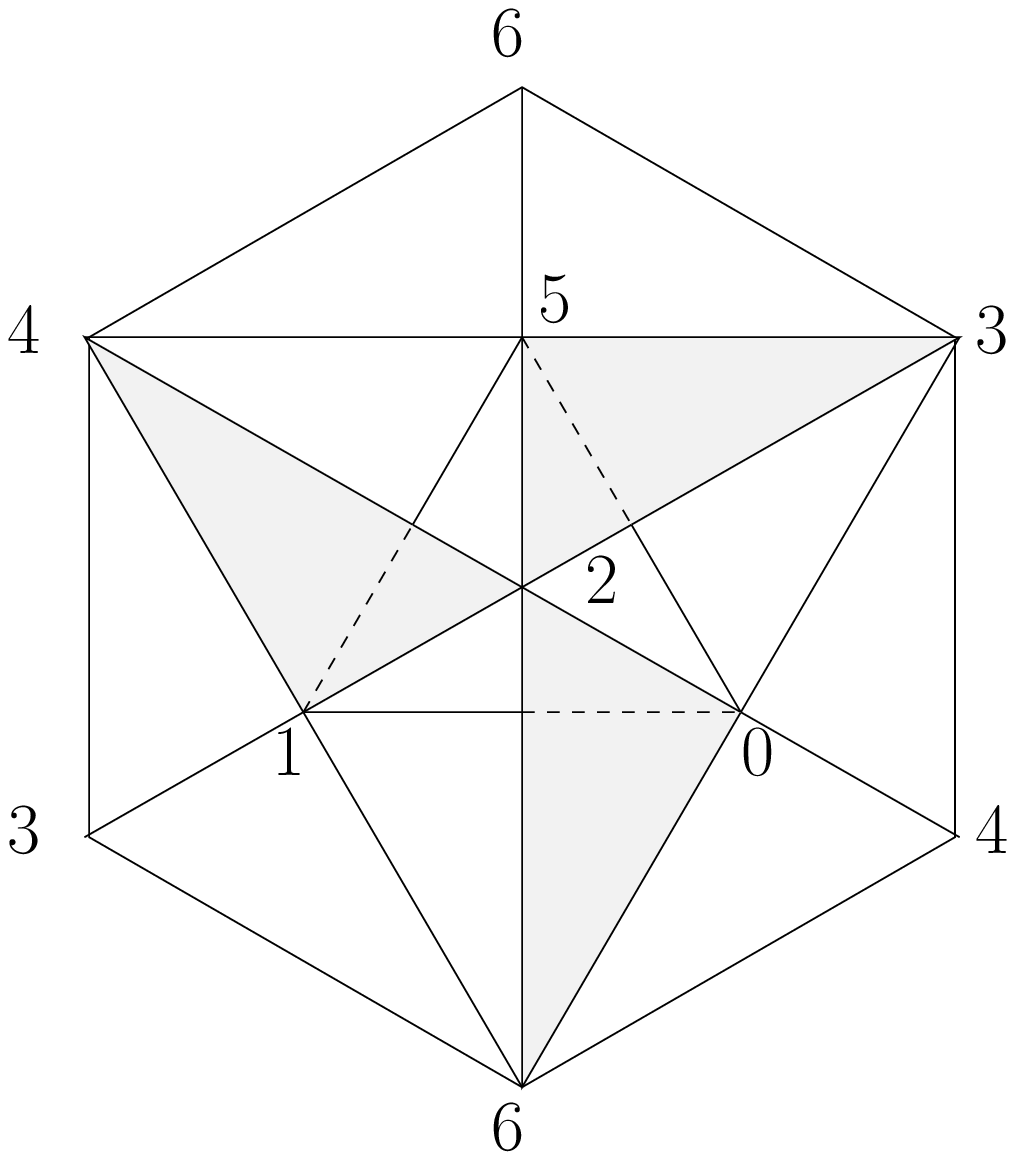}}}
  \caption{}
  \label{figure1}
\end{figure}

Let $G$ be a finite abelian group of order $n$ and let $\dn$ denote the $(n-1)$-simplex on the vertex set $G$. The \emph{sum complex} $X_{A,k}$ associated to a subset $A \subset G$ and $k < n$, is the $k$-dimensional simplicial complex obtained
by taking the full $(k-1)$-skeleton of $\dn$, together with all $(k+1)$-subsets $\sigma \subset G$ that satisfy $\sum_{x \in \sigma} x \in A$.
\\
{\bf Example:} Let $G=\Int_7$ be the cyclic group of order $7$, and let $A=\{0,1,3\}$. The sum complex $X_{A,2}$ is depicted in Figure \ref{fig:xa}). Note that $X_{A,2}$
is obtained from a $7$-point triangulation of the real projective plane $\Rea\PP^2$ (Figure \ref{fig:rp2}) by adding the faces $\{2,3,5\}$, $\{0,2,6\}$ and $\{1,2,4\}$. $X_{A,2}$ is clearly homotopy equivalent to $\Rea\PP^2$.

The sum complexes $X_{A,k}$ may be viewed as $k$-dimensional analogues of Cayley graphs over $G$. They were defined and studied (for cyclic groups) in \cite{LMR10,M14}, where some of their combinatorial and topological properties were established. For example, for $G=\Int_p$, the cyclic group of prime order $p$, the homology of $X_{A,k}$ was determined in \cite{LMR10} for coefficient fields $\FF$ of  characteristic coprime to $p$, and in \cite{M14} for general $\FF$. In particular, for $\FF=\Com$ we have the following
\begin{theorem}[\cite{LMR10,M14}]
\label{t:lmr}
Let $A \subset \Int_p$ such that $|A|=m$. Then for $1 \leq k < p$
\begin{equation*}
\label{hmodp}
\dim \tilde{H}_{k-1}(X_{A,k};\Com)=
\left\{
\begin{array}{cl}
0 & ~{\rm if}~~m \geq k+1, \\
(1-\frac{m}{k+1})\binom{p-1}{k} & ~{\rm if}~~m \leq k+1.
\end{array}
\right.~~
\end{equation*}
\end{theorem}
For a simplicial complex $X$ and $k \geq -1$ let $C^k(X)$ denote the space of complex valued
simplicial $k$-cochains of $X$ and let $d_k:C^k(X)
\rightarrow C^{k+1}(X)$ denote the coboundary operator. For $k
\geq 0$ define the reduced $k$-th Laplacian of $X$ by
$L_k(X)=d_{k-1}d_{k-1}^*+d_k^*d_k$ (see section \ref{s:hodge} for
details). The minimal eigenvalue of $L_k(X)$, denoted by $\mu_k(X)$, is the \emph{$k$-th spectral gap} of $X$.

Theorem \ref{t:lmr} implies that if $A$ is a subset of $G=\Int_p$ of size $|A|=m \geq k+1$, then $\tilde{H}_{k-1}(X_{A,k};\Com)=0$ and
hence $\mu_{k-1}(X_{A,k})>0$. Returning to the case of general finite abelian groups $G$, it is then natural to ask for better lower bounds on
the spectral gap $\mu_{k-1}(X_{A,k})$. Note that any $(k-1)$-simplex $\sigma \in \dn$
is contained in at most $m$ simplices of $X_{A,k}$ of dimension $k$, and therefore
$\mu_{k-1}(X_{A,k}) \leq m+k$ (see (\ref{e:upmu2}) in Section \ref{s:hodge}). Let $\log$ denote natural logarithm. Our main result asserts, roughly speaking, that if $k\geq 1$ and  $\epsilon>0$ are fixed and $A$ is a random subset of $G$ of size $m=\lceil c(k,\epsilon) \log n \rceil$, then
 $\mu_{k-1}(X_{A,k}) >(1-\epsilon) m$ asymptotically almost surely (a.a.s.).
The precise statement is as follows.
\begin{theorem}
\label{t:logp}
Let $k$ and $\epsilon>0$ be fixed. Let $G$ be an abelian group of order $
n> \frac{2^{10} k^8}{\epsilon^8}$, and
let $A$ be a random subset of $G$ of size
$m=\lceil\frac{4k^2\log n}{\epsilon^2}\rceil$. Then
\begin{equation*}
\label{e:sprob}
\pr\big[~\mu_{k-1}(X_{A,k}) < (1-\epsilon)m~\big] < \frac{6}{n}.
\end{equation*}
\end{theorem}

\noindent
{\bf Remarks:}
\\
1. Alon and Roichman \cite{AR94} proved that for any $\epsilon>0$ there exists a constant $c(\epsilon)>0$ such that for any group $G$ of order $n$, if $S$ is a random subset of $G$ of size
$\lceil c(\epsilon) \log n\rceil$ and $m=|S \cup S^{-1}|$, then the spectral gap of the $m$-regular Cayley graph
${\rm C}(G,S \cup S^{-1})$ is a.a.s. at least $(1-\epsilon)m$. Theorem \ref{t:logp} may be viewed as a sort of high dimensional analogue of the Alon-Roichman theorem for abelian groups.
\\
2. For $0 \leq q \leq 1$ let $Y_k(n,q)$ denote the probability space of random complexes obtained by taking the full $(k-1)$-skeleton of $\dn$ and then adding each $k$-simplex independently with probability $q$. Let $d=q(n-k)$ denote the expected number of $k$-simplices containing a fixed $(k-1)$-simplex. Gundert and Wagner \cite{GW16} proved that for any $\delta>0$ there exists a $C=C(\delta)$ such that if $q\geq \frac{(k+\delta)\log n}{n}$, then $Y \in Y_k(n,q)$ satisfies
a.a.s. $\mu_{k-1}(Y) \geq (1-\frac{C}{\sqrt{d}})d$.
\\

The paper is organized as follows. In Section \ref{s:hodge} we recall
some basic properties of high dimensional Laplacians and their eigenvalues.
In Section \ref{s:ftc} we study the Fourier images of $(k-1)$-cocycles of sum complexes,
and obtain a lower bound (Theorem \ref{t:lbmk}) on $\mu_{k-1}(X_{A,k})$, in terms of the Fourier transform of the indicator function of $A$. This bound is the key ingredient in the proof of Theorem \ref{t:logp} given in Section \ref{s:pfmain}. We conclude in Section \ref{s:con} with some remarks and open problems.

\section{Laplacians and their Eigenvalues}
\label{s:hodge}

Let $X$ be a finite simplicial complex on the vertex
set $V$. Let $X^{(k)}=\{\sigma \in X:\dim \sigma \leq k\}$ be the $k$-th skeleton of $X$, and let $X(k)$ denote the set of $k$-dimensional simplices in
$X$, each taken with an arbitrary but fixed orientation. The face numbers of $X$
are $f_k(X)=|X(k)|$. A simplicial $k$-cochain is a complex valued skew-symmetric function on
all ordered $k$-simplices of $X$. For $k \geq 0$ let $C^k(X)$
denote the space of $k$-cochains on $X$. The $i$-face of an
ordered $(k+1)$-simplex $\sigma=[v_0,\ldots,v_{k+1}]$ is the
ordered $k$-simplex
$\sigma_i=[v_0,\ldots,\widehat{v_i},\ldots,v_{k+1}]$. The
coboundary operator $d_k:C^k(X) \rightarrow C^{k+1}(X)$ is given
by $$d_k \phi (\sigma)=\sum_{i=0}^{k+1} (-1)^i \phi
(\sigma_i)~~.$$ It will be convenient to augment the cochain
complex $\{C^i(X)\}_{i=0}^{\infty}$ with the $(-1)$-degree term
$C^{-1}(X)=\Com$ with the coboundary map $d_{-1}:C^{-1}(X)
\rightarrow C^0(X)$ given by $d_{-1}(a)(v)=a$ for $a \in \Com~,~v
\in V$. Let $Z^k(X)= \ker (d_k)$ denote the space of $k$-cocycles
and let $B^k(X)={\rm Im}(d_{k-1})$ denote the space of
$k$-coboundaries. For $k \geq 0$ let $\thh^k(X)=Z^k(X)/B^k(X)~$
denote the $k$-th reduced cohomology group of $X$ with complex
coefficients. For each $k \geq -1$ endow $C^k(X)$ with the
standard inner product $(\phi,\psi)_X=\sum_{\sigma \in X(k)}
\phi(\sigma)\overline{\psi(\sigma)}~~$ and the corresponding $L^2$ norm
$||\phi||_X=(\phi,\phi)^{1/2}$.
\\ Let $d_k^*:C^{k+1}(X) \rightarrow C^k(X)$
denote the adjoint of $d_k$ with respect to these standard
inner products. The reduced $k$-th Laplacian of $X$ is the mapping
$$L_k(X)=d_{k-1}d_{k-1}^*+d_k^*d_k : C^k(X) \rightarrow
C^k(X).$$
The $k$-th Laplacian $L_k(X)$ is a positive semi-definite Hermitian operator on $C^k(X)$. Its
minimal eigenvalue, denoted by $\mu_k(X)$, is the \emph{$k$-th spectral gap} of $X$.
For two ordered simplices $\alpha \subset \beta$
let $(\beta:\alpha) \in \{\pm 1\}$ denote the incidence number between $\beta$ and $\alpha$. Let $\deg(\beta)$ denote the number of simplices $\gamma$ of dimensional $\dim \beta+1$ that contain $\beta$ . For an ordered $k$-simplex $\sigma=[v_0,\ldots,v_k] \in X(k)$, let $1_{\sigma} \in C^k(X)$ be the indicator $k$-cochain of $\sigma$, i.e. $1_{\sigma}(u_0,\ldots,u_k)={\rm sign}(\pi)$ if $u_i=v_{\pi(i)}$
for some permutation $\pi \in S_{k+1}$, and zero otherwise.
By a simple computation (see e.g. (3.4) in \cite{DR02}), the matrix representation of $L_k$ with respect to the standard basis $\{1_{\sigma}\}_{\sigma \in X(k)}$ of $C^k(X)$ is given by
\begin{equation}
\label{e:matform}
L_k(X)(\sigma,\tau)= \left\{
\begin{array}{ll}
      \deg(\sigma)+k+1 & \sigma=\tau, \\
       (\sigma :\sigma \cap \tau)\cdot (\tau:\sigma \cap \tau)
       & |\sigma \cap \tau|=k~,~\sigma \cup \tau \not\in X, \\
       0 & {\rm otherwise}.
\end{array}
\right.~~
\end{equation}
\noindent
{\bf Remarks:}
\\
(i) By (\ref{e:matform})
\begin{equation*}
\label{e:upmu1}
{\rm tr} \, L_{k}(X)=\sum_{\sigma\in X(k)} (\deg(\sigma)+k+1)=(k+2)f_{k+1}(X)+(k+1)f_k(X).
\end{equation*}
Hence
\begin{equation}
\label{e:upmu2}
\begin{split}
\mu_k(X) &\leq \frac{{\rm tr} \, L_k(X)}{f_k(X)}
\leq (k+2)\frac{f_{k+1}(X)}{f_k(X)}+k+1 \\
&\leq \max_{\sigma \in X(k)}\deg(\sigma)+k+1.
\end{split}
\end{equation}
\noindent
(ii) The matrix representation of $L_0(X)$
is equal to $J+L$, where $J$ is the $V \times V$ all ones matrix, and $L$ is the graph Laplacian of the $1$-skeleton $X^{(1)}$ of $X$. In particular, $\mu_0(X)$ is equal to the graphical spectral gap $\lambda_2(X^{(1)})$.
\ \\ \\
In the rest of this section we record some well known properties of the coboundary operators and Laplacians on the $(n-1)$-simplex $\dn$ (Claim \ref{c:simplex}), and on subcomplexes of $\dn$ that contain its full $(k-1)$-skeleton (Proposition \ref{p:ulap}). Let ${\rm I}$ denote the identity operator on $C^{k-1}(\dn)$.
\begin{claim}
\label{c:simplex}
$~$
\\
(i) The $(k-1)$-Laplacian on $\dn$ satisfies $L_{k-1}(\dn)=n\cdot {\rm I}$.  \\
(ii) There is an orthogonal decomposition
\begin{equation*}
\label{e:orth}
C^{k-1}(\dn)=\ker d_{k-2}^* \oplus \img d_{k-2}.
\end{equation*}
(iii)
The operators $P={\rm I}-\frac{1}{n}d_{k-2}d_{k-2}^*$ and $Q=\frac{1}{n}d_{k-2}d_{k-2}^*$ are, respectively, the orthogonal projections of
$C^{k-1}(\dn)$ onto $\ker d_{k-2}^*$ and onto $\img d_{k-2}$.
\end{claim}
\noindent
{\bf Proof.} Part (i) follows from (\ref{e:matform}). Next observe that
$\ker d_{k-2}^* \perp \img d_{k-2}$ and
$$\dim \ker d_{k-2}^*+\dim \img d_{k-2}=\dim \ker d_{k-2}^*+\dim \img d_{k-2}^*=\dim C^{k-1}(\dn).$$
This implies (ii).
\\
(iii)  By (i):
$$n \cdot {\rm I}=L_{k-1}(\dn)=d_{k-2}d_{k-2}^*+d_{k-1}^*d_{k-1},$$
and hence
\begin{equation*}
\label{e:sim1}
n d_{k-2}^*=d_{k-2}^*d_{k-2}d_{k-2}^*+d_{k-2}^*d_{k-1}^*d_{k-1}=d_{k-2}^*d_{k-2}d_{k-2}^*.
\end{equation*}
It follows that
$$d_{k-2}^*P=d_{k-2}^*-\frac{1}{n}d_{k-2}^*d_{k-2}d_{k-2}^*=0,$$
and therefore $\img P \subset \ker d_{k-2}^*$. Since clearly $\img Q \subset \img d_{k-2}$,
it follows that $P$ is the projection onto $\ker_{k-2}^*$ and $Q$ is the projection onto $
\img d_{k-2}$.
{\enp}
\ \\ \\
The variational characterization of the eigenvalues of Hermitian operators implies that for any complex $X$
\begin{equation}
\label{e:minev}
\begin{split}
\mu_{k-1}(X)&= \min\left\{\frac{(L_{k-1}\phi,\phi)_X}{(\phi,\phi)_X}:
0 \neq \phi \in C^{k-1}(X)\right\} \\
&=\min\left\{\frac{\|d_{k-2}^*\phi\|_X^2+\|d_{k-1}\phi\|_X^2}{\|\phi\|_X^2}:
0 \neq \phi \in C^{k-1}(X)\right\}.
\end{split}
\end{equation}
When $X$ contains the full $(k-1)$-skeleton we have the following stronger statement.
\begin{proposition}
\label{p:ulap}
Let $\dn^{(k-1)} \subset X \subset \dn$. Then
\begin{equation}
\label{e:ulap}
\mu_{k-1}(X)
=\min\left\{\frac{\|d_{k-1}\phi\|_X^2}{\|\phi\|_X^2}:
0 \neq \phi \in \ker d_{k-2}^* \right\}.
\end{equation}
\end{proposition}
\noindent
{\bf Proof.} The $\leq$ statement in (\ref{e:ulap}) follows directly from (\ref{e:minev}).
We thus have to show the reverse inequality. First note that if $\psi \in C^{k-1}(X)$ then by Claim \ref{c:simplex}(i)
\begin{equation}
\label{e:maxev}
\begin{split}
\|d_{k-1} \psi\|_X^2 &\leq \|d_{k-1} \psi\|_{\dn}^2 \\
&\leq \|d_{k-2}^* \psi\|_{\dn}^2+\|d_{k-1} \psi\|_{\dn}^2 \\
&= (L_{k-1}(\dn) \psi,\psi)_{\dn}=n \|\psi\|_X^2.
\end{split}
\end{equation}
Furthermore, if $\phi \in C^{k-1}(X)$ then by Claim \ref{c:simplex}(iii)
\begin{equation}
\label{e:dpq}
d_{k-1}\phi=d_{k-1}P\phi+d_{k-1}Q\phi=d_{k-1}P\phi,
\end{equation}
and
\begin{equation}
\label{e:dpq1}
\begin{split}
\|d_{k-2}^*\phi\|_X^2 &= (d_{k-2}^*\phi,d_{k-2}^*\phi)_X
=(\phi,d_{k-2}d_{k-2}^*\phi)_X \\
&=n(\phi,Q\phi)_X=n\|Q\phi\|_X^2.
\end{split}
\end{equation}
It follows that
\begin{equation*}
\begin{split}
\mu_{k-1}(X)&=\min\left\{\frac{\|d_{k-2}^*\phi\|_X^2+\|d_{k-1}\phi\|_X^2}{\|\phi\|_X^2}:
0 \neq \phi \in C^{k-1}(X)\right\} \\
&=\min\left\{\frac{n\|Q\phi\|_X^2+\|d_{k-1}(P\phi)\|_X^2}{\|Q\phi\|_X^2+\|P\phi\|_X^2}:
0 \neq \phi \in C^{k-1}(X)\right\} \\
&\geq \min\left\{\frac{\|d_{k-1}(P\phi)\|_X^2}{\|P\phi\|_X^2}:
0 \neq \phi \in C^{k-1}(X)\right\} \\
&= \min\left\{\frac{\|d_{k-1}\psi\|_X^2}{\|\psi\|_X^2}:
0 \neq \psi \in \ker d_{k-2}^*\right\},
\end{split}
\end{equation*}
where the first equality is (\ref{e:minev}), the second equality follows from (\ref{e:dpq}) and
(\ref{e:dpq1}), the third inequality follows from (\ref{e:maxev}) with $\psi=P\phi$, and the last equality is a consequence of Claim \ref{c:simplex}(iii).
{\enp}

\section{Fourier Transform and Spectral Gaps}
\label{s:ftc}

Let $G$ be a finite abelian group of order $n$. Let $\ccl(G)$ denote the space of complex valued functions on $G$ with the standard inner product
$(\phi,\psi)=\sum_{x \in G}\phi(x) \overline{\psi(x)}$ and the corresponding $L^2$ norm $\|\phi\|=(\phi,\phi)^{1/2}$. Let $\ghat$ be the character group of $G$.
The \emph{Fourier Transform} of $\phi \in \ccl(G)$, is the function
$\widehat{\phi} \in \ccl(\ghat)$ whose value on the character $\chi \in \ghat$ is given by
$\phat(\chi)=\sum_{x \in G} \phi(x) \chi(-x)$. For $\phi,\psi \in \ccl(G)$ we have the Parseval identity $(\phihat,\psihat)=n (\phi,\psi)$, and in particular
$\|\phihat\|^2=n\|\phi\|^2$.

Let $G^k$ denote the direct product $G \times \cdots \times G$ ($k$ times). The character group 
$\widehat{G^k}$ is naturally identified with $\ghat^k$.
Let $\tcl(G^k)$ denote the subspace of skew-symmetric functions in $\ccl(G^k)$.
Then $\tcl(G^k)$ is mapped by the Fourier transform onto $\tcl(\ghat^k)$.
Recall that $\dn$ is the simplex on the vertex set $G$, and
let $X \subset \dn$ be a simplicial complex that contains the full $(k-1)$-skeleton of $\dn$.
As sets, we will identify $C^{k-1}(X)=C^{k-1}(\dn)$ with $\tcl(G^k)$.
Note, however, that the inner products and norms defined on $C^{k-1}(\dn)$ and on $\tcl(G^k)$
differ by multiplicative constants:
If $\phi,\psi \in C^{k-1}(\dn)=\tcl(G^k)$ then
$ (\phi,\psi)=k! (\phi,\psi)_{\dn}$ and $\|\phi\|=\sqrt{k!} \|\phi\|_{\dn}$.
\ \\ \\
Let $A \subset G$ and let $k < n=|G|$. Let $\chi_0 \in \ghat$ denote the trivial character of $G$
and let $\ghat_+=\ghat \setminus \{\chi_0\}$.
Let $1_A \in \ccl(G)$ denote the indicator function of $A$, i.e. $1_A(x)=1$ if $x \in A$ and zero otherwise. Then $\widehat{1_A}(\eta)=\sum_{a \in A}\eta(-a)$ for $\eta \in \ghat$.
The main result of this section is the following lower bound on the spectral gap of $X_{A,k}$.
\begin{theorem}
\label{t:lbmk}
\begin{equation*}
\label{e:lbmk}
\mu_{k-1}(X_{A,k}) \geq |A|-k \max\{ |\widehat{1_A}(\eta)|: \eta \in \ghat_+\}.
\end{equation*}
\end{theorem}
\ \\ \\
The proof of Theorem \ref{t:lbmk} will be based on two preliminary results, Propositions \ref{p:bbb}
and \ref{p:mesti}. The first of these is the following Fourier theoretic characterization of
$\ker d_{k-2}^* \subset C^{k-1}(\dn)$.
\begin{proposition}
\label{p:bbb}
Let $\phi \in C^{k-1}(\dn)=\tcl(G^k)$. Then $d_{k-2}^*\phi =0$ iff $\supp(\phat) \subset (\ghat_+)^k$.
\end{proposition}
\noindent
{\bf Proof.} If $d_{k-2}^*\phi=0$ then for all $(x_1,\ldots,x_{k-1}) \in G^{k-1}$:
$$
0=d_{k-2}^*\phi(x_1,\ldots,x_{k-1})=\sum_{x_0 \in G} \phi(x_0,x_1,\ldots,x_{k-1}).
$$
Let $(\chi_1,\ldots,\chi_{k-1})$ be an arbitrary element of $\ghat^{k-1}$ and write
$\chi=(\chi_0,\chi_1,\ldots,\chi_{k-1}) \in \ghat^{k}$. Then
\begin{equation*}
\begin{split}
\phat(\chi) &= \sum_{(x_0,\ldots,x_{k-1}) \in G^k} \phi(x_0,x_1,\ldots,x_{k-1})
\prod_{j=1}^{k-1} \chi_j(-x_j)  \\
&=\sum_{(x_1,\ldots,x_{k-1}) \in G^{k-1}} \left(\sum_{x_0 \in G} \phi(x_0,x_1,\ldots,x_{k-1}) \right)
\prod_{j=1}^{k-1} \chi_j(-x_j) =0.
\end{split}
\end{equation*}
The skew-symmetry of $\phat$ thus implies that $\supp(\phat) \subset (\ghat_+)^k$.
The other direction is similar.
{\enp}
\ \\ \\
For the rest of this section let $X=X_{A,k}$. Fix $\phi\in C^{k-1}(X)=\tcl(G^k)$.
 Our next step is to obtain a lower bound on $\|d_{k-1}\phi\|_X$ via the Fourier transform $\widehat{d_{k-1}\phi}$. For $a \in G$ define a function $f_a \in \tcl(G^k)$ by
\begin{equation*}
\begin{split}
f_a(x_1,\ldots,x_k)&=d_{k-1}\phi\bigl(a-\sum_{i=1}^{k} x_i,x_1,\ldots,x_k\bigr) \\
&=\phi(x_1,\ldots,x_k)+ \sum_{i=1}^k (-1)^i
\phi\bigl(a-\sum_{j=1}^{k}
x_j,x_1,\ldots,\hat{x_i},\ldots,x_k\bigr).
\end{split}
\end{equation*}
By the Parseval identity
\begin{equation}
\label{e:dphi}
\begin{split}
\|d_{k-1} \phi\|_X^2 &= \sum_{\tau \in X(k)} |d_{k-1}\phi(\tau)|^2 \\
&= \frac{1}{(k+1)!} \sum_{\{(x_0,\ldots,x_k) \in G^{k+1}: \{x_0,\ldots,x_k\} \in X\}}
|d_{k-1}\phi(x_0,\ldots,x_k)|^2 \\
&= \frac{1}{(k+1)!} \sum_{a \in A} \sum_{x=(x_1,\ldots,x_k)\in G^k} |d_{k-1}\phi(a-\sum_{i=1}^k x_i,x_1,\ldots,x_k)|^2 \\
&= \frac{1}{(k+1)!}\sum_{a \in A} \sum_{x \in G^k} |f_a(x)|^2 \\
&= \frac{1}{n^k (k+1)!}\sum_{a \in A} \sum_{\chi \in \ghat^k} |\widehat{f_a}(\chi)|^2.
\end{split}
\end{equation}
We next find an expression for $\widehat{f_a}(\chi)$.
Let $T$ be the
automorphism of $\ghat^k$ given by
$$
T(\chi_1,\ldots,\chi_k)=(\chi_2 \chi_1^{-1},\ldots,\chi_{k}\chi_1^{-1},\chi_1^{-1})~.$$
Then $T^{k+1}=I$ and for $1 \leq i \leq k$
\begin{equation}
\label{e:ti}
T^i(\chi_1,\ldots,\chi_k)=(\chi_{i+1}\chi_i^{-1},\ldots,\chi_k\chi_i^{-1},\chi_i^{-1},\chi_1 \chi_i^{-1},
\ldots,\chi_{i-1}\chi_i^{-1}).
\end{equation}
\noindent
The following result is a slight extension of Claim 2.2 in \cite{LMR10}.
Recall that $\chi_0$ is the trivial character of $G$.
\begin{claim}
\label{c:ftfa}
Let $\chi=(\chi_1,\ldots,\chi_k) \in \ghatk$. Then
\begin{equation}
\label{faone}
 \widehat{f_a}(\chi)=\sum_{i=0}^{k}(-1)^{ki}
\chi_i(-a) \phat (T^{i} \chi).
\end{equation}
\end{claim}
\noindent
{\bf Proof:} For $1 \leq i \leq k$ let
$\psi_i \in \ccl(G^k)$ be given by
$$
\psi_i(x_1,\ldots,x_k)=
\phi\bigl(a-\sum_{j=1}^{k}
x_j,x_1,\ldots,\hat{x_i},\ldots,x_k\bigr).$$
Then
$$
\widehat{\psi_i}(\chi)=
\sum_{(x_1,\ldots,x_k) \in G^k}
\phi\bigl(a-\sum_{j=1}^{k}
x_j,x_1,\ldots,\hat{x_i},\ldots,x_k\bigr)
\prod_{j=1}^k \chi_j(-x_j)~~. $$
Substituting
$$
y_j=\left\{
\begin{array}{ll}
    a-\sum_{\ell=1}^k x_{\ell} & j=1, \\
    x_{j-1} & 2 \leq j \leq i, \\
    x_j & i+1 \leq j \leq k,
\end{array}
\right.
$$
it follows that
$$
\prod_{j=1}^k \chi_j(-x_j)=\chi_i^{-1}(a-y_1)\prod_{j=2}^i (\chi_i^{-1}\chi_{j-1})(-y_j) \prod_{j=i+1}^k (\chi_i^{-1}\chi_j)(-y_j)~.$$
Therefore
\begin{equation}
\label{e:psifour}
\begin{split}
\widehat{\psi_i}(\chi)&=\chi_i(-a)
\sum_{y=(y_1,\ldots,y_k)\in G^k}
\phi(y)\chi_i^{-1}(-y_1)\prod_{j=2}^i (\chi_{j-1} \chi_i^{-1}) (-y_j) \prod_{j=i+1}^k 
(\chi_j \chi_i^{-1})(-y_j) \\
&=\chi_i(-a)
\widehat{\phi}
(\chi_i^{-1},\chi_1\chi_i^{-1},\ldots,\chi_{i-1}\chi_i^{-1},\chi_{i+1}
\chi_i^{-1},\ldots,\chi_k\chi_i^{-1}) \\
&=\chi_i(-a) (-1)^{i(k-i)} \widehat{\phi}(T^i \chi)~.
\end{split}
\end{equation}
Now (\ref{faone}) follows from (\ref{e:psifour}) since
$f_a=\phi+\sum_{i=1}^k (-1)^i \psi_i$.
{\enp}

\noindent
For $\phi \in \tcl(G^k)$ and $\chi=(\chi_1,\ldots,\chi_k) \in \ghat^k$ let
\begin{equation*}
\label{e:dchi}
\begin{split}
&D(\phi,\chi)= \{\chi_i \chi_j^{-1}: 0 \leq i < j \leq k~,
~\phat(T^i \chi)\phat(T^j\chi)\neq 0\} \\
&=\{\chi_j^{-1}:1 \leq j \leq k~,~\phat(\chi)\phat(T^j\chi) \neq 0\} \cup
\{\chi_i \chi_j^{-1}: 1 \leq i < j \leq k~,
~\phat(T^i \chi)\phat(T^j\chi)\neq 0\}.
\end{split}
\end{equation*}
Let $D(\phi)=\bigcup_{\chi \in \ghat^k} D(\phi,\chi) \subset \ghat$.
The main ingredient in the proof of Theorem \ref{t:lbmk} is the following
\begin{proposition}
\label{p:mesti}
$$
\|d_{k-1} \phi\|_X^2 \geq \left(|A|- k\max_{\eta \in D(\phi)} |\widehat{1_A}(\eta)|\right)\|\phi\|_X^2.
$$
\end{proposition}
\noindent
{\bf Proof.}
Let $\chi=(\chi_1,\ldots,\chi_k) \in \ghat^k$.
By Claim \ref{c:ftfa} 
\begin{equation}
\label{suma}
\begin{split}
&\sum_{a \in A} |\widehat{f_a}(\chi)|^2=\sum_{a \in A} |\sum_{i=0}^{k}(-1)^{ki}
\chi_i(-a) \phat (T^{i} \chi)|^2 \\
&=\sum_{a \in A} \sum_{i,j=0}^{k}(-1)^{k(i+j)}
(\chi_i \chi_j^{-1})(-a) \phat (T^{i} \chi)\overline{\phat (T^{j} \chi)} \\
&= |A|\sum_{i=0}^k  |\phat (T^{i} \chi)|^2+ 2\, \text{Re}\sum_{a \in A} \sum_{0 \leq i<j \leq k}(-1)^{k(i+j)}
(\chi_i \chi_j^{-1})(-a)\phat(T^i\chi)\overline{\phat(T^j\chi)} \\
&= |A|\sum_{i=0}^k  |\phat (T^{i} \chi)|^2+ 2\, \text{Re}\sum_{0 \leq i<j \leq k}(-1)^{k(i+j)}
\widehat{1_A}(\chi_i\chi_j^{-1})\phat(T^i \chi)\overline{\phat(T^j\chi)} \\
&\geq |A|\sum_{i=0}^k  |\phat (T^{i} \chi)|^2 -2\max_{\eta \in D(\phi,\chi)} |\widehat{1_A}(\eta)|\sum_{0 \leq i<j \leq k} |\phat(T^i \chi)|\cdot|\phat(T^j \chi)|.
\end{split}
\end{equation}
Using (\ref{e:dphi}) and summing (\ref{suma}) over all $\chi \in \ghat^k$ it follows that
\begin{equation*}
\label{sumau}
\begin{split}
&n^k (k+1)!\|d_{k-1} \phi\|_X^2 = \sum_{a \in A} \sum_{\chi \in \ghat^k} |\widehat{f_a}(\chi)|^2 \\
&\geq  |A|\sum_{i=0}^k  \sum_{\chi \in \ghat^k} |\phat (T^{i} \chi)|^2 -2\max_{\eta \in D(\phi)} |\widehat{1_A}(\eta)|\sum_{0 \leq i<j \leq k}
\sum_{\chi \in \ghat^k} |\phat(T^i \chi)|\cdot|\phat(T^j\chi)| \\
&\geq (k+1)|A|\sum_{\chi \in \ghat^k} |\phat(\chi)|^2- k(k+1)\max_{\eta \in D(\phi)} |\widehat{1_A}(\eta)| \, \sum_{\chi \in \ghat^k} |\phat(\chi)|^2 \\
&=(k+1)n^k \sum_{x \in G^k} |\phi(x)|^2 \left(|A|-k\max_{\eta \in D(\phi)} |\widehat{1_A}(\eta)|\right) \\
&= (k+1)!n^k\|\phi\|_X^2\left(|A|- k\max_{\eta \in D(\phi)} |\widehat{1_A}(\eta)|\right).
\end{split}
\end{equation*}
{\enp}
\noindent
{\bf Proof of Theorem \ref{t:lbmk}.} Let $0 \neq \phi \in C^{k-1}(X_{A,k})=\tcl(G^k)$ such that $d_{k-2}^*\phi=0$. Proposition \ref{p:bbb} implies that
\begin{equation}
\label{e:supht}
\supp(\phat) \subset (\ghat_+)^k.
\end{equation}
We claim that
$\chi_0 \not\in D(\phi)$. Suppose to the contrary that $\chi_0 \in D(\phi)$. Then there exists
a $\chi=(\chi_1,\ldots,\chi_k) \in \ghat^k$ such that $\chi_0 \in D(\phi,\chi)$,
i.e. $\chi_i=\chi_j$ for some
$0 \leq i < j \leq k$ such that $\phat(T^i \chi)\phat(T^j \chi) \neq 0$. We consider two cases:
\begin{itemize}
\item
If $i=0$ then $\chi_j=\chi_i=\chi_0$ and therefore
$$0 \neq \phat(T^i\chi)=\phat(\chi)=
\phat(\chi_1,\ldots,\chi_{j-1},\chi_0,\chi_{j+1},\ldots,\chi_k),$$
in contradiction of (\ref{e:supht}).
\item
If $i \geq 1$ then $\chi_j \chi_i^{-1}=\chi_0$, and by (\ref{e:ti})
\begin{equation*}
\label{e:e:supht1}
\begin{split}
0&\neq \phat(T^i \chi)=
\phat(\chi_{i+1}\chi_i^{-1},\ldots,\chi_k\chi_i^{-1},\chi_i^{-1},\chi_1 \chi_i^{-1},
\ldots,\chi_{i-1}\chi_i^{-1}) \\
& \phat (\chi_{i+1}\chi_i^{-1},\ldots,\chi_{j-1} \chi_i^{-1},\chi_0,\chi_{j+1}\chi_i^{-1},\ldots,
\chi_k\chi_i^{-1},\chi_i^{-1},\chi_1 \chi_i^{-1},
\ldots,\chi_{i-1}\chi_i^{-1}),
\end{split}
\end{equation*}
again in contradiction of (\ref{e:supht}).
\end{itemize}
We have thus shown that $D(\phi) \subset \ghat_+$.
Combining Propositions \ref{p:ulap} and \ref{p:mesti} we obtain
\begin{equation*}
\label{e:fin}
\begin{split}
\mu_{k-1}&(X_{A,k})
=\min\left\{\frac{\|d_{k-1}\phi\|_X^2}{\|\phi\|_X^2}:
0 \neq \phi \in \ker d_{k-2}^* \right\} \\
&\geq \min\left\{
\frac{\left(|A|- k\max_{\eta \in D(\phi)}
|\widehat{1_A}(\eta)|\right)\|\phi\|_X^2 }{\|\phi\|_X^2}:
0 \neq \phi \in \ker d_{k-2}^* \right\}  \\
&\geq  |A|- k\max_{\eta \in \ghat_+} |\widehat{1_A}(\eta)|.
\end{split}
\end{equation*}
{\enp}

\section{Proof of Theorem \ref{t:logp}}
\label{s:pfmain}

Let $k \geq 1$ and $0<\epsilon<1$ be fixed, and let $n> \frac{2^{10} k^8}{\epsilon^8}$,
$m=\lceil\frac{4k^2\log n}{\epsilon^2}\rceil$. Let $G$ be an abelian group of order $n$ and
let $\Omega$ denote the uniform probability space of all $m$-subsets of $G$.
Suppose that $A \in \Omega$ satisfies $|\widehat{1_A}(\eta)| \leq \frac{\epsilon m}{k}$ for all
$\eta \in \ghat_+$.  Then by Theorem \ref{t:lbmk}
\begin{equation*}
\begin{split}
\mu_{k-1}(X_{A,k}) &\geq |A|-k \max\{ |\widehat{1_A}(\eta)|: \eta \in \ghat_+\} \\
& \geq |A|-k \cdot \frac{\epsilon m}{k}=(1-\epsilon)m.
\end{split}
\end{equation*}
Theorem \ref{t:logp} will therefore follow from
\begin{proposition}
\label{p:cher}
\begin{equation*}
\label{e:cher}
\pr_{\Omega}\left[~A \in \Omega~:~\max_{\eta \in \ghat_+}
|\widehat{1_A}(\eta)|>\frac{\epsilon m}{k}~\right]< \frac{6}{n}.
\end{equation*}
\end{proposition}
\noindent
{\bf Proof.}
Let $\eta \in \ghat_+$ be fixed and let $\lambda=\frac{\epsilon m}{k}$.
Let $\Omega'$ denote the uniform probability space $G^m$, and for $1 \leq i  \leq m$ let $X_i$ be the random variable defined on $\omega'=(a_1,\ldots,a_m) \in \Omega'$ by
$X_i(\omega')=\eta(-a_i)$. The $X_i$'s are of course independent and satisfy $|X_i|=1$. Furthermore, as $\eta \neq \chi_0$, the expectation of $X_i$ satisfies  $E_{\Omega'}[X_i]=\frac{1}{n}
\sum_{x \in G} \eta(-x)=0$. Hence, by the Chernoff bound (see e.g. Theorem A.1.16 in \cite{AS08})
\begin{equation}
\label{e:chernoff1}
\begin{split}
&\pr_{\Omega'}\left[~\omega' \in \Omega'~:~|\sum_{i=1}^m X_i(\omega')|>\lambda~\right]
< 2\exp\left(-\frac{\lambda^2}{2m}\right) \\
&=2\exp\left(-\frac{\epsilon^2m}{2k^2}\right) \leq \frac{2}{n^2}.
\end{split}
\end{equation}
Let $\Omega''=\{(a_1,\ldots,a_m) \in G^m: a_i \neq a_j \text{~for~} i \neq j \}$
denote the subspace of $\Omega'$ consisting of all sequences in $G^m$ with pairwise distinct elements. Note that the assumption $n > 2^{10} k^8 \epsilon^{-8}$ implies that
\begin{equation}
\label{e:mp}
\frac{m^2}{n-m} <1.
\end{equation}
Combining (\ref{e:chernoff1}) and (\ref{e:mp}) we obtain
\begin{equation}
\label{e:chernoff2}
\begin{split}
&\pr_{\Omega}\left[A \in \Omega:|\widehat{1_A}(\eta)|>\frac{\epsilon m}{k}~\right] \\
&=\pr_{\Omega''}\left[~\omega'' \in \Omega'':|\sum_{i=1}^m X_i(\omega'')|
>\frac{\epsilon m}{k}~\right] \\
&\leq \pr_{\Omega'}\left[~\omega' \in \Omega'~:~|\sum_{i=1}^m X_i(\omega')| >
\frac{\epsilon m}{k}~\right]
\cdot \left(\pr_{\Omega'}[~\Omega''~]\right)^{-1} \\
&< \frac{2}{n^2} \cdot \prod_{i=1}^m \frac{n}{n-i+1} \leq \frac{2}{n^2} \cdot
\left(\frac{n}{n-m}\right)^m \\
&\leq \frac{2}{n^2} \cdot \exp\left(\frac{m^2}{n-m}\right) < \frac{6}{n^2}.
\end{split}
\end{equation}
Using the union bound together with (\ref{e:chernoff2}) it follows that
$$
\pr_{\Omega}\left[~\max_{\eta \in \ghat_+}|\widehat{1_A}(\eta)|>\frac{\epsilon m}{k}~\right]< \frac{6}{n}.
$$
{\enp}

\section{Concluding Remarks}
\label{s:con}
In this paper we studied the $(k-1)$-spectral gap of sum complexes $X_{A,k}$ over a finite abelian group $G$. Our main results include a Fourier theoretic lower bound on $\mu_{k-1}(X_{A,k})$,
and a proof that if $A$ is a random subset of $G$ size of $O(\log |G|)$, then $X_{A,k}$ has a nearly optimal $(k-1)$-th spectral gap.
Our work suggests some more questions regarding sum complexes:

\begin{itemize}
\item
Theorem \ref{t:logp} implies that if $G$ is an abelian group of order $n$, then
$G$ contains many subsets $A$ of size $m=O_{k,\epsilon}(\log n)$ such that
$\mu_{k-1}(X_{A,k})\geq (1-\epsilon)m$. As is often the case with probabilistic existence proofs, it would be interesting to give explicit constructions for such $A$'s.
For $G=\Int_2^{\ell}$, such a construction
follows from the work of Alon and Roichman. Indeed, they observed (see Proposition 4 in \cite{AR94}) that by results of \cite{ABNNR92},
there is an absolute constant $c>0$ such that for any $\epsilon>0$ and $\ell$, there is an explicitly constructed $A_{\ell} \subset \Int_2^{\ell}$ of size
$m\leq \frac{c k^3 \ell}{\epsilon^3}=\frac{c k^3 \log_2 |G|}{\epsilon^3}$, such that
$$
|\widehat{1_{A_\ell}}(v)|=|\sum_{a \in A}(-1)^{a \cdot v}| \leq \frac{\epsilon m}{k}
$$
for all $0 \neq v \in \Int_2^{\ell}$. Theorem \ref{t:lbmk} then implies that
$\mu_{k-1}(X_{A_{\ell},k}) \geq (1-\epsilon)m$.
\\
It would be interesting to find explicit constructions with $|A|=O(\log |G|)$ for other groups $G$ as well, in particular for the cyclic group $\Int_p$.

\item
Consider the following non-abelian version of sum complexes. Let $G$ be a finite group of order $n$
and let $A \subset G$. For $1 \leq i \leq k+1$ let $V_i$ be the $0$-dimensional complex on the set $G \times \{i\}$, and let $T_{n,k}$ be the join $V_1*\cdots*V_{k+1}$.
The complex $R_{A,k}$ is obtained by taking the $(k-1)$-skeleton of $T_{n,k}$, together with all
$k$-simplices $\sigma=\{(x_1,1),\ldots,(x_{k+1},k+1)\} \in T_{n,k}$ such that
$x_1\cdots x_{k+1} \in A$. One may ask whether there is an analogue of Theorem \ref{t:logp}
for the complexes $R_{A,k}$, i.e. is there a constant $c_1(k,\epsilon)>0$ such that
if $A$ is a random subset of $G$ of size $m=\lceil c_1(k,\epsilon) \log n \rceil$, then
a.a.s. $\mu_{k-1}(R_{A,k}) >(1-\epsilon) m$.
\end{itemize}

\end{document}